\documentclass{amsart}
\usepackage{amsmath,amsthm,amssymb,bm}
\usepackage{color}
\usepackage{mathtools}
\mathtoolsset{showonlyrefs=true} 

\usepackage{tikz}
\usetikzlibrary{intersections,patterns}

\usepackage{hyperref}

\numberwithin{equation}{section}

\newtheorem{definition}{Definition}[section]

\newtheorem{theorem}[definition]{Theorem}
\newtheorem{lemma}[definition]{Lemma}

\newtheorem{remark}[definition]{Remark}

\newcommand{\re}{\mathbb R}

\DeclareMathOperator{\diver}{div}
\DeclareMathOperator{\rot}{rot}
\DeclareMathOperator{\supp}{supp}

\title[A remark on Liouville-type theorems for stationary NS]{
A remark on Liouville-type theorems for the stationary Navier--Stokes equations in $\mathbb{R}^3$
under pointwise conditions on the head pressure
} 
\author{Hideo Kozono, Yutaka Terasawa and Yuta Wakasugi}

\address[H. Kozono]{Department of Mathematics, Faculty of Science and Engineering,
Waseda University, Tokyo 169--8555, Japan, 
Mathematical Research Center for Co-creative Society, Tohoku University, 
Sendai 980-8578, Japan}
\email[H. Kozono]{kozono@waseda.jp, hideokozono@tohoku.ac.jp}
\address[Y. Terasawa]{Graduate School of Mathematics, Nagoya University,
Furocho Chikusaku Nagoya 464-8602, Japan}
\email[Y. Terasawa]{yutaka@math.nagoya-u.ac.jp}
\address[Y. Wakasugi]{Graduate School of Advanced Science and Engineering,
Hiroshima University,
Higashi-Hiroshima, 739-8527, Japan}
\email[Y. Wakasugi]{wakasugi@hiroshima-u.ac.jp}


\dedicatory{Dedicated to Professor Yoshikazu Giga on the award of the Kodaira Kunihiko Prize}

\begin{document}
\maketitle

\begin{abstract}
We study the stationary Navier--Stokes equations in the three-dimensional whole space
and the Liouville-type problem under the condition of finite Dirichlet integral. 
We give a remark on the recent result by Chae \cite[Theorem 1.2]{Ch26},
in which an upper bound of the gradient of the head pressure and
a lower bound of the head pressure with certain rates were assumed 
to prove triviality of solutions.
In this paper, we extend the range of admissible rates by two differenct arguments.
One relies on the pointwise decay of head pressure and the classical argument of Galdi \cite{Ga}.
The other is similar to the argument of Chae \cite{Ch20, Ch25} and is based on
the maximum principle for the head pressure,
the integral identity between the vorticity and the gradient of the head pressure on the upper level set of the head pressure,
and the coarea formula.
This also gives a more refined proof for \cite[Theorem 1.2]{Ch26}.
\end{abstract}

\section{Introduction}
In this paper, we consider the stationary Navier--Stokes equations in the whole space $\mathbb{R}^3$:
\begin{equation}\label{eq:NS}
	\left\{ \begin{alignedat}{3}
	&-\Delta v + (v \cdot \nabla) v + \nabla p = 0,\\
	&\diver v = 0 
	\end{alignedat} \right.
	\quad \text{in} \ \mathbb{R}^3.
\end{equation}
Here,
$v = v(x) = (v_1(x), v_2(x), v_3(x))$
and
$p=p(x)$
are the velocity and the pressure at the point
$x = (x_1, x_2, x_3) \in \mathbb{R}^3$,
respectively.
We also impose on $v$ the condition of finite Dirichlet integral 
\begin{equation}\label{eq:Diri}
	\int_{\mathbb{R}^3} | \nabla v(x) |^2 \,dx < + \infty
\end{equation}
with the uniform decay at the spatial infinity
\begin{equation}\label{eq:decay}
	\lim _{|x| \to \infty}|v(x)|=0 \quad \text { uniformly in } x.
\end{equation}

In Galdi \cite[Remark X.9.4, p.729]{Ga} and Tsai \cite[Conjecture 2.5]{Ts},
the following Liouville-type problem was stated:
\textit{Are there nonidentically vanishing smooth solutions $v, p$ to the Navier--Stokes problem
\eqref{eq:NS} satisfying \eqref{eq:Diri} and \eqref{eq:decay} ?}

For this Liouville-type problem, there are many partial results
by imposing on $v$ additional integrability or decay conditions. 
Galdi \cite{Ga} proved that if $v$ additionally satisfies
$v \in L^{9/2}(\mathbb{R}^3)$,
then $v \equiv 0$.
Chae and Wolf \cite{ChWo} gave a logarithmic improvement to Galdi's result.
Chae \cite{Ch14} showed that
$\Delta v \in L^{6/5}(\mathbb{R}^3)$
and \eqref{eq:decay} lead to triviality of $v$.
In our previous paper \cite{KoTeWa17},
it was proved that
if \eqref{eq:decay} holds and
$\limsup_{|x|\to\infty} |x|^{5/3} \rot v$
is sufficiently small,
or
$\| v \|_{L^{9/2,\infty}}$
is sufficiently small,
then $v$ must vanish identically.
After that, Seregin and Wang \cite{SeWa} gave a generalization of
$L^{9/2,\infty}$-condition to the smallness of the Morrey-type quantity
$\liminf_{R\to \infty} R^{2/3-3/q} \| v \|_{L^{q, \ell}(R/2\le |x| \le R)}$
with some
$q > 3, 3 \le \ell \le \infty$ (or $q = \ell = 3$).
They also proved that if the quantity
$\liminf_{R\to \infty} R^{\gamma-3/q} \| v \|_{L^{q, \ell}(R/2\le |x| \le R)}$
is finite for some
$12/5 < q < 3, 1 \le \ell \le \infty, \gamma > 1/3 + 1/q$,
then
$v$ must vanish identically.
Tsai \cite{Ts21} improved them by including the borderline cases
$\gamma = 1/3 + 1/q$ and $q = 12/5$
for the case of Lebesgue spaces, that is, $\ell = q$.
Cho, Neustupa and Yang \cite{ChoNeYa} further improved the range of
$q$ to
$3/2<q<3$.

On the other hand, Seregin \cite{Se} obtained a sufficient condition
$v \in BMO^{-1}(\mathbb{R}^3) \cap L^6(\mathbb{R}^3)$
implying
$v \equiv 0$.
Chamorro, Jarr\'in and Lemari\'e-Rieusset \cite{ChaJaLe}
generalized it in the sense that
they replaced
$BMO^{-1}(\mathbb{R}^3)$
by a Morrey-type space with the same scaling
and $L^6$ by similar Morrey-type spaces with more general scaling.
Recently, Tan \cite{Ta} proved that if $v$ satisfies \eqref{eq:Diri} and
$v \in \dot{B}^{-1}_{\infty, \infty} (\mathbb{R}^3)$,
then $v$ must be zero.
We note that the space
$\dot{B}^{-1}_{\infty, \infty} (\mathbb{R}^3)$
has the same scaling as
$BMO^{-1}(\mathbb{R}^3)$.
The proof of \cite{Ta} is based on the paraproduct decomposition and the analysis of
the low frequency part of $v$.
Several sufficient conditions regarding
the low frequency part of
$v$
in the Fourier space was also studied by Lerner \cite{Ler} recently.
\par
\bigskip
In this paper, we focus on the results and methods by Chae \cite{Ch20, Ch21, Ch25, Ch26},
in which additional conditions in terms of the head pressure
\[
	Q(x) = \frac{1}{2}|v(x)|^2 + p(x)
\]
had been imposed.
In \cite{Ch20}, the following sufficient condition (and its logarithmic improvements)
was introduced:
\begin{equation}\label{eq:Ch20}
	\int_{\mathbb{R}^3} \frac{|\nabla Q|^2}{|Q|} \,dx < + \infty.
\end{equation}
Here, as we will see later, we may assume that
$Q < 0$
everywhere without loss of generality.
In \cite{Ch21}, it was shown that if either
\begin{equation}\label{eq:Ch21}
	\sup_{x\in\mathbb{R}^3} \frac{|v(x)|^2}{|Q(x)|} < + \infty
	\quad \text{or} \quad
	\sup_{x\in\mathbb{R}^3} \frac{|p(x)|}{|Q(x)|} < + \infty,
\end{equation}
holds, then $v \equiv 0$
(see also Bang--Yang \cite{BaYa} for the logarithmic improvement
for the denominator).
In \cite{Ch25}, this result was generalized by replacing
$|Q(x)|$
in the denominator by
$\xi(|Q(x)|)$
with a function
$\xi(s)$
satisfying the Osgood condition
\[
	\int_0^1 \frac{1}{\xi(s)} \,ds = + \infty.
\]
Finally, in \cite{Ch26}, the following two types of pointwise assumptions were introduced:
\begin{itemize}
\item[(i)]
There exist
$\alpha, \beta >0$
with
$\beta \ge \alpha/2$
and positive constants
$C, R$
such that
\begin{equation}
	|Q(x)| \geq \frac{C\|Q\|_{L^{\infty}}}{|x|^\alpha} \quad (\forall |x|>R)
	\quad \text{and} \quad
	|v(x)|=O\left(\frac{1}{|x|^\beta}\right) \ \text { as } \ |x| \rightarrow+\infty ; 
\end{equation}
\item[(ii)]
There exist
$\alpha, \beta >0$
with
$\beta \ge 2\alpha$
and positive constants $C, R$ such that
\begin{equation}
	|Q(x)| \geq \frac{C\|Q\|_{L^{\infty}}}{|x|^\alpha} \quad (\forall |x|>R)
	\quad \text{and} \quad
	|\nabla Q(x)|=O\left(\frac{1}{|x|^\beta}\right) \ \text { as } \ |x| \rightarrow+\infty.
\end{equation}
\end{itemize}
It was proved that if either (i) or (ii) is satisfied as well as the condition of finite Dirichlet integral \eqref{eq:Diri}, then
$v$
must be a constant vector
(in particular, $v \equiv 0$ under (i)).

\begin{remark}\label{rem:Q}
Roughly speaking, the above results by \cite{Ch20, Ch21, Ch25, Ch26} require some estimates for $|Q(x)|$ from below.
The reason is that, for example, if
$|Q(x)| \ge C|x|^{-\alpha}$ with some $\alpha>0$,
then the sub-level set $\{ x ; |Q(x)| < \lambda \}$ with small $\lambda>0$
is included in the set
$\{ x ; |x| > C \lambda^{-\frac{1}{\alpha}} \}$.
\end{remark}
Our main result concerns the condition (ii) of \cite{Ch26} above,
and extends the range of
$(\alpha, \beta)$
by the argument of Galdi \cite{Ga} and Chae \cite{Ch20, Ch25}. 
\par
\bigskip
Before stating the main result, we prepare for some properties of $Q(x)$.
By \cite[Theorem X.5.1]{Ga},
for a smooth solution $(v,p)$ to \eqref{eq:NS} with the condition  \eqref{eq:Diri} of finite Dirichlet integral,
there exist a constant vector
$v_{\infty} \in \mathbb{R}^3$
and a constant
$p_1 \in \mathbb{R}$
such that 
\begin{equation}\label{eq:vp:infty}
	\lim_{|x|\to \infty} \left( |v(x) - v_{\infty}| + |p(x)-p_1| \right) = 0
	\quad \text{uniformly in} \ x.
\end{equation}
\begin{remark}\label{rem:v:p}
By \cite[Theorem X.7.2]{Ga} with \eqref{eq:NS}, \eqref{eq:Diri} and \eqref{eq:vp:infty},
it follows that if $v_{\infty} \neq 0$, then $v \equiv v_{\infty}$ and
the Liouville-type theorem is proved.
Therefore, in what follows, we assume that $v_{\infty} = 0$,
that is,
$\lim_{|x| \to \infty} v (x) = 0$
uniformly in $x$.
Moreover, since the equation \eqref{eq:NS} does not change even if we replace
$p$ by $p-p_1$,
we may assume that
$p$ satisfies
$\lim_{|x|\to \infty} p(x) = 0$
uniformly in $x$.
We also remark that, by the condition  \eqref{eq:Diri} of finite Dirichlet integral  and the Sobolev embedding, we have
$v \in L^6(\mathbb{R}^3)$.
\end{remark}
From this remark, we may assume that $Q$ satisfies
\begin{equation}\label{eq:Q:decay}
	\lim_{|x|\to \infty} | Q(x) | = 0 \quad \text{uniformly in} \ x.
\end{equation}

\begin{theorem}\label{thm}
Let
$(v, p)$ be a smooth solution to \eqref{eq:NS} with the condition \eqref{eq:Diri} of finite Dirichlet integral,
and let
$Q$
be the associate head pressure satisfying \eqref{eq:Q:decay}.

{\rm (i)}
Assume that there exists
$\beta \ge 5/2$
such that
\[
	|\nabla Q(x)|=O\left(\frac{1}{|x|^\beta}\right) \ \text { as } \ |x| \rightarrow+\infty.
\]
Then, $v \equiv 0$.
\par
{\rm (ii)} 
Assume that there exist
$\alpha, \beta > 0$
with
\[
	\beta \ge 2\alpha
	\quad \text{or} \quad
	\beta > \frac{\alpha+3}{2}
\]
and positive constants
$C, R$
such that
\begin{equation}
	|Q(x)| \geq \frac{C\|Q\|_{L^{\infty}}}{|x|^\alpha} \quad (\forall |x|>R)
	\quad \text{and} \quad
	|\nabla Q(x)|=O\left(\frac{1}{|x|^\beta}\right) \ \text { as } \ |x| \rightarrow+\infty.
\end{equation}
 Then, $v \equiv 0$.
\end{theorem}
\begin{remark}
Chae \cite[Theorem 1.2]{Ch26} proved the same conclusion under the condition
$\beta \ge 2\alpha$.
Therefore, the above theorem extends the admissible range of $(\alpha, \beta)$.
In Figure 1, we draw the three regions of
$(\alpha, \beta)$
corresponding to the results of Chae \cite{Ch26}, Theorem \ref{thm} (i) and (ii).
\end{remark}

\begin{figure}[htbp]
\centering

\begin{tikzpicture}

\draw (0,2.5)--(4,2.5) node[right] {$\beta = \frac{5}{2}$}; 
\fill[lightgray] (0,2.5)--(4,2.5)--(4,4)--(0,4)--cycle;
\node[left] at (0,2.5) {$\frac{5}{2}$};

\draw (0,0)--(2,4) node[above]{$\beta=2\alpha$}; 
\fill[pattern=horizontal lines, pattern color=gray] (0,0)--(2,4)--(0,4)--cycle; 

\draw[densely dashed] (0,1.5)--(4,3.5) node[right] {$\beta = \frac{\alpha+3}{2}$}; 
\fill[pattern=vertical lines, pattern color=gray] (0,1.5)--(4,3.5)--(4,4)--(0,4)--cycle;
\node[left] at (0,1.5) {$\frac{3}{2}$};

\node[above left] at (3,3) {{(\textit{ii})}}; 
\node[above] at (3.5,2.5) {{(\textit{i})}};
\node[right] at (0, 3.5) {Chae \cite{Ch26}};

\draw[->, thick] (-0.5,0)--(4,0) node[right] {$\alpha$}; 
\draw[->, thick] (0,-0.5)--(0,4) node[left] {$\beta$};

\end{tikzpicture}

\caption{}
\end{figure}

\section{Proof of the main theorem}
It can be seen that the head pressure
$Q(x)$
obeys the equation
\begin{equation}\label{eq:Q:diffeq}
	-\Delta Q + v \cdot \nabla Q = - |\omega|^2,
\end{equation}
where
$\omega = \rot v$.
One can verify this by a straightforward computation with noting that
\begin{equation}\label{eq:Lap:p}
	-\Delta p = \sum_{j,k=1}^3 \partial_j v_k \partial_k v_j,
\end{equation}
which follows from taking the divergence to the equation \eqref{eq:NS}.

We first note that $v \equiv 0$ is equivalent to $Q \equiv 0$.
Indeed,
if $v \equiv 0$, then by the equation \eqref{eq:NS} we obtain $\nabla p = 0$, namely, $p$ is a constant.
Moreover, by Remark \ref{rem:v:p} (or \eqref{eq:Q:decay}), we have $\lim_{|x| \to \infty} p(x)=0$,
which implies $p \equiv 0$, and hence, $Q \equiv 0$.
Conversely, if $Q \equiv 0$, then from \eqref{eq:Q:diffeq} we obtain $\omega \equiv 0$.
Hence, the formula
$\rot (\rot v) = \nabla (\diver v) - \Delta v$
implies that
$v$
is a harmonic function.
Thus, the condition \eqref{eq:vp:infty} and the classical Liouville theorem for the harmonic function lead to
$v \equiv 0$.

Therefore, to prove Theorem \ref{thm}, it suffices to show $Q \equiv 0$.
We will prove it by contradiction, and in what follows, we suppose
$Q \not\equiv 0$.
Then, since the right-hand side of \eqref{eq:Q:diffeq} is nonpositive,
$Q$
satisfies the maximum principle and does not achieve the maximum in
$\mathbb{R}^3$.
This and the behavior  \eqref{eq:Q:decay} at the infinity imply
$Q \le 0$.
Moreover, applying the maximum principle again, we deduce that
$Q < 0$
in $\mathbb{R}^3$.
\par
(i) We first treat the case
$\beta \ge 5/2$.
By the assumption
$\nabla Q(x) = O(|x|^{-\beta})$ with $\beta > 0$, we see that\footnote{
From this estimate we remark that if $\beta > \alpha+1$, then it immediately contradicts
the assumption $|Q(x)| \ge C \| Q \|_{L^{\infty}} |x|^{-\alpha}$,
while we do not need this fact in the subsequent proof.
}
\[
	|Q(x)| \le \int_{|x|}^{\infty} \left| \nabla Q \left( r \frac{x}{|x|} \right) \right| \,dr
	\le C |x|^{-\beta + 1}. 
\]
In particular, since
$\beta \ge 5/2$,
we have
$|Q(x)| \le C |x|^{-3/2}$.
Then, we apply the argument of Galdi \cite[Remark X.9.4]{Ga}.
Let
$\psi \in C_0^{\infty}(\mathbb{R}^3)$
be such that
\[
	0\le \psi(x) \le 1, \quad
	\psi(x)=
	\begin{cases} 1 &(|x| \le 1),\\ 0 &(|x| \ge 2) \end{cases}
\]
and let
\[
	\psi_R(x) = \psi \left( \frac{x}{R} \right)
\]
with a parameter
$R>0$.
We remark that
\begin{equation}\label{eq:psi:est}
	|\nabla \psi_R (x) | \le C R^{-1},\quad
	|\Delta \psi_R(x) | \le C R^{-2}
\end{equation}
and
\begin{equation}\label{eq:psi:supp}
	\supp \nabla \psi_R, \supp \Delta \psi_R \subset B_{2R} \setminus B_R
\end{equation}
hold, where
$B_R$
denotes the ball centered at the origin with radius $R$.
Taking the inner product of the equation \eqref{eq:NS} and
$v \psi_R$,
integrating it over $\mathbb{R}^3$,
and using integration by parts,
we obtain the identity
\begin{align}
	\int_{\mathbb{R}^3} |\nabla v|^2 \psi_R \,dx
	&= \frac{1}{2} \int_{\mathbb{R}^3} |v|^2 \Delta \psi_R \,dx
	+ \int_{\mathbb{R}^3} Q v \cdot \nabla \psi_R \,dx \\
	&=: I_1 + I_2.
\end{align}
By the fact $v \in L^6(\mathbb{R}^3)$, the properties \eqref{eq:psi:est}, \eqref{eq:psi:supp}
and the H\"{o}lder inequality, we easily obtain
\begin{align}
	|I_1| &\le C \| v \|_{L^6(B_{2R}\setminus B_R)}^2 \| \Delta \psi_R \|_{L^{3/2}(B_{2R}\setminus B_R)} \\
	&\le C \| v \|_{L^6(B_{2R}\setminus B_R)}^2.
\end{align}
Similarly, the estimate
$|Q(x)| \le C |x|^{-3/2}$
implies
\begin{align}
	|I_2| & \le C \| Q \|_{L^2(B_{2R}\setminus B_R)} \| v \|_{L^6(B_{2R}\setminus B_R)} 
			\| \nabla \psi_R \|_{L^3(B_{2R}\setminus B_R)} \\
		&\le C \| v \|_{L^6(B_{2R}\setminus B_R)}.
\end{align}
Therefore, we conclude
\[
	\int_{\mathbb{R}^3} |\nabla v|^2 \psi_R \,dx
	\le C \left(  \| v \|_{L^6(B_{2R}\setminus B_R)}^2 +  \| v \|_{L^6(B_{2R}\setminus B_R)} \right).
\]
By letting $R \to \infty$,
the left-hand side tends to
$\| \nabla v \|_{L^2(\mathbb{R}^3)}^2$,
while the right-hand side tends to $0$.
Thus, we have
$\nabla v \equiv 0$, and $v$ must be a constant vector.
Since $v$ satisfies
$\lim_{|x|\to \infty} v(x) = 0$,
we conclude
$v \equiv 0$.
From this and the fact that $v \equiv 0$ is equivalent to $Q \equiv 0$ proved in the beginning of this section,
we obtain $Q \equiv 0$.
This contradicts the assumption $Q \not\equiv 0$ and completes the proof.
\par
\vspace{1mm}
(ii) Next, we consider the other case
$\beta \ge 2\alpha$ or $\beta > \frac{\alpha+3}{2}$.
Although the case $\beta \ge 2\alpha$ was already proved by Chae \cite{Ch26},
we give here a more refined proof based on the argument of Chae \cite{Ch20, Ch25}.
We start with the following important identity:
\begin{lemma}\label{lem:Q:omega}
Let
$\omega = \rot v$.
For almost every
$\lambda \in (0, \| Q \|_{L^{\infty}})$,
we have
\begin{equation}\label{eq:Q:omega}
	\int_{\{ |Q(x)|>\lambda \}} |\omega|^2 \,dx
	= \int_{ \{ |Q(x)|=\lambda \}} |\nabla Q| \,dS.
\end{equation}
\end{lemma}
\begin{proof}
By the Morse--Sard theorem
(see e.g., \cite{Mi} for a simple proof),
wee see that almost every
$\lambda \in (0, \| Q \|_{L^{\infty}})$
is a regular value of $|Q| = -Q$.
For such a value, we see 
by (\ref{eq:Q:decay}) that 
the set
$\{ x ; |Q(x)| = \lambda\}$
is a union of smooth closed surfaces so that we can apply the divergence theorem in the upper level set
$\{ x ; |Q(x)| > \lambda \}$.
Integrating \eqref{eq:Q:diffeq} over this upper level set, we obtain
\begin{align*}
	\int_{ \{ |Q(x)| > \lambda \}} |\omega|^2 \,dx
	&=
	\int_{ \{ |Q(x)| > \lambda \}} \Delta Q \,dx
	- \int_{ \{ |Q(x)| > \lambda \}} v \cdot \nabla Q \,dx \\
	&=
	\int_{ \{ |Q(x)| = \lambda \}} \frac{\partial Q}{\partial \nu} \,dS
	- \int_{ \{ |Q(x)| = \lambda \}} v \cdot \nu Q \,dS,
\end{align*}
where
$\nu$
is the outer unit normal on the surface
$\{ |Q(x)| = \lambda \}$,
and we have used the divergence free condition
$\diver v = 0$
to the second term.
Applying the divergence theorem and
$\diver v = 0$
again,
we see that the second term is
\[
	- \int_{ \{ |Q(x)| = \lambda \}} v \cdot \nu Q \,dS
	= \lambda \int_{ \{ |Q(x)| = \lambda \}} v \cdot \nu \,dS
	= \lambda \int_{ \{ |Q(x)| > \lambda \}} \diver v \,dx = 0.
\]
Since
$\nu = \nabla Q/|\nabla Q|$, the first term is rewritten as 
\[
	\int_{ \{ |Q(x)| = \lambda \}} \frac{\partial Q}{\partial \nu} \,dS
	= \int_{ \{ |Q(x)| = \lambda \}} |\nabla Q| \,dS.
\]
This completes the proof of Lemma \ref{lem:Q:omega}.
\end{proof}
Let us replace the notation
$\lambda$ in Lemma \ref{lem:Q:omega} by $\mu$.
Multiplying the identity \eqref{eq:Q:omega} by
$1/ \mu$
and then integrating it over
$(\lambda/2, \lambda)$
with respect to $\mu$,
we have
\[
	\int_{\lambda/2}^{\lambda} \frac{1}{\mu} \int_{\{ |Q(x)| > \mu\} } |\omega|^2 \,dx d\mu
	=
	\int_{\lambda/2}^{\lambda} \frac{1}{\mu} \int_{\{ |Q(x)| = \mu\} } |\nabla Q| \,dS d\mu.
\]
When $\lambda \to 0$,
the left-hand side obviously satisfies
\begin{equation}\label{eqn:identity}
	\lim_{\lambda \to 0} \int_{\lambda/2}^{\lambda} \frac{1}{\mu} \int_{\{ |Q(x)| > \mu\} } |\omega|^2 \,dx d\mu
	= \log 2 \cdot \| \omega \|_{L^2(\mathbb{R}^3)}^2.
\end{equation}
For the right-hand side, we apply the coarea formula
(see e.g., \cite{EvGa})
to obtain
\begin{align*}
	\int_{\lambda/2}^{\lambda} \frac{1}{\mu} \int_{\{ |Q(x)| = \mu\} } |\nabla Q| \,dS d\mu
	&=
	\int_{\lambda/2}^{\lambda}  \int_{\{ |Q(x)| = \mu\} } \frac{|\nabla Q|}{|Q|} \,dS d\mu \\
	&=
	\int_{\{ \frac{\lambda}{2} < |Q(x)| < \lambda \} } \frac{|\nabla Q|^2}{|Q|} \,dx.
\end{align*}

Now, we divide the cases.
First we consider the case
$\beta > \frac{\alpha+3}{2}$.
Recall the assumptions
$\nabla Q = O(|x|^{-\beta})$ and $|Q| \ge C \| Q \|_{L^{\infty}} |x|^{-\alpha}$
for $|x| > R$.
By (\ref{eq:Q:decay}) we note that the set
$\{ x ; \frac{\lambda}{2} < |Q(x)| < \lambda \}$
is included in
$\{ x ; |x| > R \}$
for sufficiently small $\lambda$.
Therefore, it holds that 
\[
	\int_{\{ \frac{\lambda}{2} < |Q(x)| < \lambda \} } \frac{|\nabla Q|^2}{|Q|} \,dx
	\le
	C \int_{c\lambda^{-1/\alpha}}^{\infty} r^{\alpha-2\beta +2} \,dr,
\]
where we have also used the fact that
$C |x|^{-\alpha} \le |Q(x)| < \lambda$ implies $|x| \ge c \lambda^{-1/\alpha}$.
Since $\alpha -2 \beta + 2 < -1$, implied by $\beta > \frac{\alpha + 3}{2}$, 
 the right-hand side tends to $0$ as $\lambda \to 0$.
Therefore, we conclude from \eqref{eqn:identity} that 
$\omega \equiv 0$,
which gives
$v \equiv 0$.
Since $v \equiv 0$ is equivalent to $Q \equiv 0$, we obtain $Q \equiv 0$.
This contradicts the assumption $Q \not\equiv 0$ and completes the proof.

\par
Next, we consider the case $\beta \ge 2\alpha$. 
Let us first show that 
$\nabla Q \in L^{3/2}(\mathbb{R}^3)$.  
From the representation of $Q$, we have that  
$\nabla Q = \frac{1}{2} \nabla |v|^2 + \nabla p$.
The H\"{o}lder inequality immediately gives
\[
	\| \nabla |v|^2 \|_{L^{3/2}}
	\le C \| \nabla v \|_{L^2} \| v \|_{L^6}. 
\]
Moreover, from \eqref{eq:Lap:p}, we have 
\[
	\nabla p = \nabla \sum_{j,k=1}^3 R_j R_k (v_j v_k),
\]
where
$R_j = \partial_j(-\Delta)^{-1/2}$
is the Riesz transform.
Since the Riesz transform is bounded on
$L^q(\mathbb{R}^3)$ for $q \in (1,\infty)$,
we deduce
\[
	\| \nabla p \|_{L^{3/2}} \le C \sum_{j,k=1}^3 \| \nabla (v_j v_k) \|_{L^{3/2}} 
	\le C \| \nabla v \|_{L^2} \| v \|_{L^6}. 
\]
This proves
$\nabla Q \in L^{3/2}(\mathbb{R}^3)$. 
\par
Therefore, since 
$\nabla Q(x) = O(|x|^{-\beta})$ and $|Q(x)| \ge C \|Q \|_{L^{\infty}} |x|^{-\alpha}$
for $|x| > R$,  and since $-\beta/2 + \alpha \le 0$, implied by $\beta \ge 2\alpha$, 
we have for small $\lambda$ that 
\begin{align*}
	\int_{\{ \frac{\lambda}{2} < |Q(x)| < \lambda \} } \frac{|\nabla Q|^2}{|Q|} \,dx
	&\le
	C \int_{\{ \frac{\lambda}{2} < |Q(x)| < \lambda \} } |x|^{-\beta/2+\alpha} |\nabla Q|^{3/2} \,dx \\
	&\le
	C \int_{\{ \frac{\lambda}{2} < |Q(x)| < \lambda \} } |\nabla Q|^{3/2} \,dx,
\end{align*}
where we have again used the fact
$\{ x ; \frac{\lambda}{2} < |Q(x)| < \lambda \} \subset \{ x ; |x| > R \}$
for sufficiently small $\lambda$. 
Since $\nabla Q \in L^{3/2}(\re^3)$, 
the right-hand side tends to $0$ as $\lambda \to 0$.
Hence again by (\ref{eqn:identity}) we see that this implies
$\omega \equiv 0$,
which yields that 
$v \equiv 0$.
This again implies $Q \equiv 0$, and it contradicts the assumption $Q \not\equiv 0$
and completes the proof of Theorem \ref{thm}.

\section*{Acknowledgement}
This work was supported by
JST Moonshot R\&D Grant Number JPMJMS25A5
and JSPS Grant-in-Aid for Scientific Research (A)-26H02004,
Research (A)-JP22H00097, Research (C)-JP24K06811
and Grant-in-Aid for Transformative Research Areas (B)-25H01453.  

 \section*{Data Availability Statement}
Data sharing is not applicable to this article as no datasets were generated or analysed during the current study.

\section*{Competing Interests}
The authors declare that they have no competing interests.


\section*{Declaration of generative AI use}
During the preparation of this manuscript,
the authors used Claude (Anthropic) to assist in identifying potentially relevant references.
All references included in the manuscript were independently verified and selected by the authors.

\end{document}